\documentclass[12pt]{article}

\usepackage{amsthm}

\usepackage{latexsym, amsmath, amssymb, epsfig}
\usepackage{graphicx}
\usepackage{epstopdf}

\DeclareGraphicsExtensions{.jpg,.pdf,.eps}

\usepackage[margin=3cm]{geometry}

\numberwithin{equation}{section}

\newtheorem{theorem}{Theorem}[section]

\newtheorem{lemma}[theorem]{Lemma}

\newcommand{\p}{\partial}

\newcommand{\eqnref}[1]{(\ref {#1})}

\def\Be{{\bf e}}

\def\Bu{{\bf u}}

\newcommand{\Gvf}{\varphi}
\newcommand{\Gg}{\gamma}

\newcommand{\Gl}{\lambda}

\newcommand{\Gm}{\mu}

\newcommand{\Gt}{\theta}

\newcommand{\Gr}{\rho}

\newcommand{\Go}{\omega}

\newcommand{\Gy}{\psi}

\newcommand{\GT}{\Theta}

\newcommand{\GO}{\Omega}

\newcommand{\beq}{\begin{equation}}
\newcommand{\eeq}{\end{equation}}

\allowdisplaybreaks
\title{Stationary viscoelastic wave fields generated by scalar wave functions}
\author{Junyong Eom$^1$, Gen Nakamura$^{2}$
\\{\small $^1$Department of Mathematics, Tohoku University, Sendai 980-8578, Japan}\quad\
\\{\small E-mail: eom.junyong.r2@dc.tohoku.ac.jp}
\\{\small $^2$Department of Mathematics, Hokkaido University, Sapporo 060-0810, Japan} \quad
\\{\small E-mail: gnaka@math.sci.hokudai.ac.jp}
}

\date{}
\begin{document}
\maketitle

\begin{abstract}
The usual Helmholtz decomposition gives a decomposition of any vector valued function into a sum of gradient of a scalar function and rotation of a vector valued function under some mild condition. In this paper we show that the vector valued function of the second term i.e. the divergence free part of this decomposition can be further decomposed into a sum of a vector valued function polarized in one component and the rotation of a vector valued function also polarized in the same component. Hence the divergence free part only depends on two scalar functions. Further we show the so called completeness of representation associated to this decomposition for the stationary wave field of a homogeneous, isotropic viscoelastic medium. That is by applying this decomposition to this wave field, we can show that each of these three scalar functions satisfies a Helmholtz equation. Our completeness of representation is useful for solving boundary value problem in a cylindrical domain for several partial differential equations of systems in mathematical physics such as stationary isotropic homogeneous elastic/viscoelastic equations of system and stationary isotropic homogeneous Maxwell equations of system. As an example, by using this completeness of representation, we give the solution formula for torsional deformation of a pendulum of cylindrical shaped homogeneous isotropic viscoelastic medium.  
\end{abstract}

\section{Introduction}\label{Sec:Introduction}${}$

The Helmholtz decomposition has played an important role to solve boundary value problems in a rectangle paralellpiped, sphere and cylinder for several partial differential equations of systems in mathematical physics such as stationary isotropic homogeneous elastic/viscoelastic equations of system and stationary isotropic Maxwell equations of system.

The aim of this paper is to provide a special Helmholtz decomposition for solutions of an isotropic homogeneous elastic/viscoelastic equations of system in a finite cylindrical domain. 
The special Helmholtz decomposition enables to represent solutions of this homogeneous viscoelastic system in a finite cylindrical shaped domain by scalar wave functions via a special Helmholtz decomposition associated to the cylindrical coordinates. This is called the {\it completeness of representation}. More precisely let $\Omega$ be a homogeneous isotropic viscoelastic medium in ${\mathbb R}^3$ with complexified Lam\'e moduli $\lambda,\,\mu$ such that their real parts $\lambda_R,\,\mu_R$ and imaginary parts $\lambda_I,\,\mu_I$ satisfy the strong convexity condition:
\begin{equation}\label{strconvexity}
\mu_R,\,\mu_I>0,\,\, 3\lambda_R+2\mu_R,\,3\lambda_I+2\mu_I>0
\end{equation}
(see \cite{R-W}).
The shape of $\Omega$ could be either a ball or cylinder. Define the stationary homogeneous isotropic viscoelastic operator $L_{\lambda,\mu}$ with density $\rho>0$ and angular frequency $\omega>0$ by
\begin{equation} \label{viscosystem}
L_{\lambda,\mu}u:=(\lambda+\mu)\text{\rm grad}\,\text{\rm div} u+\mu\Delta u+\rho\omega^2 u.
\end{equation}
Then the completeness of representation is given as follows. That is for $u$ satisfying $L_{\lambda,\mu} u=0$ in $\Omega$, there exists scalar functions $\varphi,\,\psi,\,\chi$ such that $u$ has a special Helmholtz decomposition given as
\begin{equation} \label{repsol}
u=\nabla\varphi+\nabla\times(\psi e)+\nabla\times\nabla\times(\chi e)
\end{equation}
and $\varphi,\,\psi,\,\chi$ satisfy the following Helmoholtz type equations:
\begin{equation}
\left\{
\begin{array}{l}
(\Delta+\frac{\rho\omega^2}{\lambda+2\mu})\varphi=0,\\
(\Delta+\frac{\rho\omega^2}{\mu})\psi=0,\\
(\Delta+\frac{\rho\omega^2}{\mu})\chi=0,
\end{array}
\right.
\end{equation}
where $e$ is a unit vector in the axial direction and direction of $z$ axis for the case
$\Omega$ is a ball and cylinder, respectively. If there is a global orthogonal curvilinear coordinates $(x^1,x^2,x^3)$ in a neighborhood of $\overline\Omega$, then there is a question about the completeness of representation when $e$ is replaced by $w(x^3)e_3$ with unit vector $e_3$ say in the $x^3$ direction and some auxilary function $w(x^3)$ for the seperation of variables. 

This completeness of representation was rigorously proved in \cite{chadwick} when a medium has a shape of ball. After that in \cite{Eringen-book}, a more general discussion was made for the completeness of representation. It included the case that the medium has a shape of cylinder. However it did not give any rigorous proof for this case. To the best of our knowledge, it seems that the proof for this case has been left opened for more than 40 years.

In this paper our main aim is to give a rigorous proof for the completeness of representation for the case $\Omega$ is a finite cylinder.
Hence in the rest of the paper, $\Omega$ denotes an open finite cylinder.

Needless to say that this completeness of representation \eqnref{repsol} is very useful for giving a solution formula to solve given boundary value problem. As an example, we apply the torque boundary condition(see \eqnref{traction}) and find the solution which satisfies our representation \eqnref{repsol}. In fact we used the solution formula to identified the complex shared modulus of the sample using a digitized commercial base pendulum-type viscoelastic spectrometer (PVS) (see \cite{EKNW} for its details).

The rest of this paper is organized as follows. In Section \ref{Sec:Spl-Helm-Decompo}, we give a special type of Helmholtz decomposition. Then we give the completeness of representation for this decomposition in Section \ref{Sec:Comp-Rep}. The last section is devoted to an application of the completeness of representation. 
That is we apply it to solve a boundary value problem for a cantilever homogeneous isotropic viscoelastic cylinder with a given torque at its bottom by using the completeness of representation. 

\section{Special Helmholtz Decomposition} \label{Sec:Spl-Helm-Decompo} 
For $l,m \in {\Bbb Z_+}:= {\Bbb N}\cup\{0\}$, $l \geq m$ and $p \geq 2$, we define a function space
$$  v \in W_{p, loc}^{l,m} ({\Bbb R}^3)  \Longleftrightarrow \sum_{|\alpha| + j \leq l ,  j \leq m} \int_{{\Bbb R}^3}
 |  \p_{x'}^{\alpha} \p_z^j (\eta v) (x',z) |^p dx'dz < \infty  $$
for all $\eta \in C_0^{\infty} ({\Bbb R}^3)$, where we used $x=(x,y,z)\in{\Bbb R}^3$, $x'=(x,y)$. We also define $W_{p}^{l,m} (\Omega)$ by deleting $\eta$ and replacing ${\Bbb R}^3$ by $\Omega$ in the above definition. For $p=2$,  we use the notations $H_{loc}^{l,m}({\Bbb R}^3)$, $H^{l,m}(\Omega)$. Also, we can infer that
$$  \p_{x'}^{\beta} \p_z^{k} v \in W_{p,loc}^{l - |\beta| -k, \min\{m-k,l-|\beta|-k\}} ({\Bbb R}^3) $$
which follows from 
$$ \left\{
\begin{array}{lll}
|\alpha|+|\beta|+j+k \leq l, \\
j+k \leq m. \\
\end{array}
\right.                    \Longleftrightarrow 
\left\{
\begin{array}{lll}
|\alpha|+ j \leq  l - |\beta| - k, \\
j \leq m -k. \\
\end{array}
\right.
$$ 
\begin{theorem}\label{special:helm} Assume that $u\in W_p^6({\Bbb R}^3)$, $p \geq 2$ and it has a compact support. Then $u$ admits the representation 
\begin{eqnarray}\label{SHD}
u= \nabla \varphi + \nabla \times (\psi e_z) + \nabla \times \nabla \times (\chi e_z),
\end{eqnarray}
where $\varphi\in W^7_{p, loc}({\Bbb R}^3)$, $ \psi \in W_{p, loc}^{7,4}(({\Bbb R}^3)$, $\chi \in W_{p, loc}^{8,5}(({\Bbb R}^3)$. Note that the regularities of $\psi,\,\chi$
given above imply $\nabla\times(\psi e_z),\,\,\nabla\times\nabla\times(\chi e_z)\in W_{p, loc}^{6,4}({\Bbb R}^3)$.
\end{theorem}
\noindent
{\bf Proof.}
First of all, it is well known that for any given $u \in W_p^6({\Bbb R}^3)$ with compact support, we can have the following Helmholtz decomposition \cite{Dautray}:
\begin{equation}\label{hDecomposition}
u=\nabla\varphi+\nabla\times A\,\,\text{in}\,\,{\Bbb R}^3\,\,\text{with}
\end{equation}
where $\varphi$ and $A$ are given by
\begin{equation}\label{varphi and A}
\begin{array}{l}
\varphi(x)=-(4\pi)^{-1}\int\frac{\nabla_y \cdot u(y)}{|y-x|} dy\in W^7_{p, loc}({\Bbb R}^3),\\
A(x)=(4\pi)^{-1}\int\frac{\nabla_y \times u(y)}{|y-x|} dy\in H^7_{p, loc}({\Bbb R}^3).
\end{array}
\end{equation}
where the integration are take over ${\Bbb R}^3$. 

Next in order to decompose $\nabla \times A$ into two potentials, we prepare the 
solvability of the equation for $\chi$ which satisfies 
\begin{equation}\label{equation for psi}
\Delta'\chi=-e_z\cdot (\nabla \times A) := f \in W^6_{p, loc}({\Bbb R}^3) \,\,\text{in}\,\,{\Bbb R}^2,
\end{equation}
where $\Delta'$ is the Laplacian with respect to $x'=(x_1,x_2)=(x,y)$.
Since the inner product $(A(x),b)$ of $A(x)$ with any constant vector $b$ is given as
\begin{equation}
\begin{array}{ll}
(\nabla \times A(x),b)&=(4\pi)^{-1}\int(\nabla_y \times u(y), \nabla_x \times (|y-x|^{-1}b)\,dy\\
&=-(4\pi)^{-1}\int\big(u(y), \nabla_y \times \nabla_x \times (|y-x|^{-1}b)\big)\,dy,
\end{array}
\end{equation}
we have $\partial_{x'}^\alpha (\nabla \times A(x))=O((1+|x|)^{-3-|\alpha|})\,\,(|x|\rightarrow\infty)$ for $|\alpha|\le1$.
Hence fixing $z\in{\Bbb R}$, we have 
\begin{equation}\label{estimate of A}
\partial_{x'}^\alpha (\nabla \times A(x',z))=O((1+|x'|)^{-3-|\alpha|})\,\,(|x'|\rightarrow\infty)
\end{equation}
for $|\alpha|\le1$. In the rest of this proof we will omit $(|x|\rightarrow\infty)$ and $(|x'|\rightarrow\infty)$ for simplicity.
\par
Now we will apply the following result in \cite{M} for the case $n=2$ in which we are interested. That is, if
$p, \, \delta$ satisfies the condition
\begin{equation}\label{condition for indices}
p>2, -\frac{2}{p} -1<\delta< -\frac{2}{p} \qquad \text{except}\, \, \delta + \frac{2}{p} \, \text{or} -(\delta + \frac{2}{p}) \in{\Bbb Z}_+,
\end{equation}
then $\Delta': M_{p,\delta}^2\rightarrow L_{\delta+2}^2$ is surjective and the kernel of $\Delta'$ is ${\Bbb C}$. Here $L_{\delta+2}^p = L_{\delta+2}^p ({\Bbb R}^2) $ and $M_{p,\delta}^2$ is the completion of $C_0^{\infty}({\Bbb R}^2)$ in the norm
\begin{equation}
\sum_{|\alpha|\le 2}\Vert (1+|x'|^2)^{(\delta+|\alpha|)/2}\partial_{x'}^\alpha h\Vert_{L^p({\Bbb R}^2)}
\end{equation}
with $x=(x',z)=(x,y,z)$, $\partial_{x'}^{\alpha} = \partial_{x_1}^i \partial_{x_2}^j$, $ i+j = \alpha$ for nonnegative $i,j$ and $L_{\delta+2}^p$ is the $L^2({\Bbb R}^2)$ with weight $(1+|x'|^2)^{(\delta+2)/2}$. Let $ \delta$ satisfy \eqref{condition for indices}, then we have $e_z\cdot (\nabla \times A(x)) \in L_{\delta+2}^p$ because
\begin{eqnarray*}
O\{(1+|x'|)^{\delta+2}(e_z \cdot (\nabla \times A(x',z)))\} = O((1+|x'|)^{-1+\delta}) \\
\end{eqnarray*}
and $p(-1 + \delta)<-2$ since $\delta<-2/p$. Therefore, the solution $\chi\in M_{p,\delta}^2$ of \eqref{equation for psi} exists.
\par
Let $\dot{M}_{p,\delta}=\{\zeta\in M_{p,\delta}^2: \int_{|x'|\le R}\zeta(x')\,dx'= 0\}$ with
a large enough fixed $R>0$. Then, since the kernel of $\Delta'$ is ${\Bbb C}$, 
$$
\Delta': \dot{M}_{p,\delta}\rightarrow L_{\delta+2}^p
$$ 
is a bi-continuous isomorphism by
Banach's closed graph theorem and open mapping theorem. Since the mapping
$$
{\Bbb R}\ni z\rightarrow e_z\cdot (\nabla \times A(\cdot, z)) \in L_{\delta+2}^p
$$
is locally continuous, the solution $\chi=\chi(\cdot,z)\in \dot{M}_{p,\delta}$ depends
continously on $z\in{\Bbb R}$. 
\par
Furthermore, we can show $\p_{x'} \chi \in W_{p,loc}^{7,5}({\Bbb R}^3)$. 
By using the fundamental solution for $\Delta^{'}$, we have following representation \cite{Nirenberg}, that
\begin{eqnarray*} \p_{x'} \chi (x',z) = \frac{1}{2 \pi} \int_{{\Bbb R}^2} \log (|x'-y'|) \p_{y'} f(y',z) dy'.
\end{eqnarray*}
Above representation holds because $\p_{x'} f(x',z) \in L_{\delta +1}^p ({\Bbb R}^3) \subset L^p ({\Bbb R}^3)$ when we take  $\delta+1>0$. So, since $\p_{x'} f(x',z) \in  W^5_{p, loc}({\Bbb R}^3)$, we have $\p_{x'} \chi \in W_{p,loc}^{7,5}({\Bbb R}^3)$. Therefore, we have $\chi \in W_{p,loc}^{8,5}({\Bbb R}^3)$.
\par
In a similar way, we can easily check that the existence of a solution $\psi\in \dot{M}_{2,\delta}^2$ of 
\begin{equation}
\Delta'\psi=-e_z\cdot(\nabla\times \nabla \times A(x)) := g \in W^5_{p, loc}({\Bbb R}^3)\,\,\text{in}\,\,{\Bbb R}^2
\end{equation}
and its properties which implies $\p_{x'} \psi \in W_{p,loc}^{6,4}({\Bbb R}^3)$, where $p,\delta$ satisfy \eqnref{condition for indices} and $\delta +1 >0$. So, we have $\psi \in W_{p,loc}^{7,4}({\Bbb R}^3)$.
\par
Now define $V$ by 
\begin{equation}\label{def of V}
V=\nabla \times A-\nabla\times\nabla\times(\chi e_z)-\nabla \times (\psi e_z).
\end{equation}
Then, we have the following lemma.
\begin{lemma}\label{Laplace eq for V}
\begin{equation}\label{eq for V}
(\Delta' V)(x',z)=0\,\,\text{\rm in}\,\,{\Bbb R}_{x'}^2\,\,\text{\rm for any $z\in{\Bbb R}$}.
\end{equation}
\end{lemma}
\noindent
{\bf Proof.}
Observe that for $V =  V_1 e_x + V_2 e_y + V_3 e_z$, we have 
\begin{equation}\label{several equations for V}
\left\{
\begin{array}{l}
\text{div}\, V =0,\\
e_z\cdot V=e_z\cdot (\nabla \times A)-e_z\cdot\big(\nabla\times\nabla\times(\chi e_z)\big)=e_z\cdot(\nabla \times A)+\Delta'\chi=0,\\
e_z\cdot\big(\nabla\times V\big)=e_z\cdot\big(\nabla \times \nabla\times A\big)-e_z\cdot\big(\nabla\times\nabla\times(\psi e_z)\big)=e_z\cdot\big(\nabla\times \nabla\times A\big)+\Delta'\psi=0.
\end{array}
\right.
\end{equation}
Then the proof can be finished by just observing that the first and third equations give the Cauchy Riemann system of equations for $V-(e_z\cdot V)e_z = V_1 e_x + V_2 e_y$. 
This completes the proof of Lemma\ref{Laplace eq for V}.

\bigskip
To proceed further, we recall the following result in \cite{M}. That is if $\delta'$ satisfies the condition
\begin{equation}\label{condition for indices 2}
-\frac{2}{p}<\delta'< -\frac{2}{p} +1 \qquad\text{except}\,\,\delta' + \frac{2}{p} \, \,\text{or}\,-(\delta' + \frac{2}{p}) \in {\Bbb Z}_+,
\end{equation}
then $\Delta': M_{p,\delta'}^2\rightarrow L_{\delta'+2}^p$ is injective. We will consider $\delta'=\delta+1$ for $\delta$ satisfying \eqnref{condition for indices}.
\par
Now, we show that $V \in M_{2,\delta'}^{2}$ for each fixed $z\in{\Bbb R}$. Since $\nabla \times A = O((1+|x'|)^{-3})$, we need $p(\delta' - 3)< -2$. It holds from the condition $\delta'<0$ . Hence $\nabla \times A \in  M_{2,\delta'}^{2}$.

\medskip
To examine the second term of $V$, we need a following lemma in \cite{M}.
\begin{lemma}\label{L,M relation} If $u \in L_{\delta'}^{p} $, $\Delta' u \in L_{\delta'+ 2}^{p}$, then $u \in M_{p,\delta'}^{2}$.
\end{lemma} 
For any constant vector $b$, we have from \eqref{equation for psi} that
\begin{eqnarray*} \Delta' (b \cdot \nabla \times \nabla \times (\chi e_z)) = b \cdot \nabla \times \nabla \times ( (-e_z \cdot \nabla \times A)e_z ) = O((1+|x'|)^{-5}).
\end{eqnarray*}
We can easily check that $p(-5+ \delta' + 2) < -2$. Hence $\Delta' (b \cdot \nabla \times \nabla \times (\chi e_z)) \in L_{\delta'+2}^{p}$. Also, for the potential $\chi$ of \eqnref{equation for psi}, we have $\chi \in M_{p,\delta}^2$ and by the definition of $M_{p,\delta}^2$, we have $b \cdot \nabla \times \nabla \times (\chi e_z) \in L_{\delta+2}^p \subset L_{\delta'}^{2}$ with the $\delta$ of \eqnref{condition for indices}. Therefore by Lemma \ref{L,M relation}, we have $ \nabla \times \nabla \times (\chi e_z) \in M_{p,\delta'}^{2}$.
\par
As we did in the second term, for the third term of $V$, we have
\begin{eqnarray*} \Delta' (b \cdot \nabla \times (\psi e_z)) = b \cdot \nabla \times ( (-e_z \cdot \nabla \times \nabla \times A)e_z ) = O((1+|x'|)^{-5})
\end{eqnarray*}
and hence $\Delta' (b \cdot \nabla  \times (\psi e_z)) \in L_{\delta'+2}^{p}$. Also, from $\psi \in M_{p,\delta}^2$ with $\delta$ of \eqnref{condition for indices}, we have $\nabla  \times (\psi e_z) \in L_{\delta'}^{p}$. Hence by Lemma 2.2, we have $ \nabla  \times (\psi e_z) \in M_{p,\delta'}^{2}$.
\par
Therefore we have shown $V=\nabla \times A-\nabla\times\nabla\times(\psi e_z)-\nabla \times (\chi e_z) \in M_{p,\delta'}^{2}$ which together with Lemma \ref{Laplace eq for V} immediately gives  $V=0$. This completes the proof of Theorem \ref{special:helm}.

\section{Completeness of representation} \label{Sec:Comp-Rep} 

In this section we will show the completeness of representation based on the special Helmholtz decomposition. From now on for our convenience, we will use the convention embedding $L^p$- spaces into $L^2$- spaces for $p>2$ and bounded domains. For example, $W_p^7 (\Omega) \subset H^7(\Omega)$, $W_p^{7,5} (\Omega) \subset H^{7,5}(\Omega)$. If we don't use this convention, it will be stated.
\

We first preare a Boggio's type result as lemma for further arguments.
\begin{lemma} \label{inhomogeneous Boggio 1} Let $\varphi \in H^{m+2}(\Omega), \, m \geq 0$ be a solution for $\nabla^2 ( \Box_1 \varphi  )= 0$ in $\Omega$ where $\Box_1 = \nabla^2 +c_1^2$, $c$ is constant. Then $\varphi$ can be decomposed into  $\varphi = \varphi_1 + \varphi_2$ with $\varphi_1,\varphi_2 \in H^m(\Omega)$ such that $\nabla^2 \varphi_1=0$ and $\Box_1 \varphi_2 =0$ in $\Omega$. 
\end{lemma} 
\noindent
{\bf Proof.} 
Let $\varphi = \varphi_1 + (\varphi - \varphi_1)$ such that $\varphi_1 = \frac{1}{c_1^2} \Box_1 \varphi  \in H^m(\Omega)$. 
Then we can easily check that $\nabla^2 \varphi_1=0$ and 
$\Box_1(\varphi - \varphi_1) =  \Box_1( - \frac{1}{c_1^2} \nabla^2 \varphi ) = 0 $ in $\Omega$. 
Hence we can finish the  proof by setting $\varphi_2 = \varphi - \varphi_1 \in H^m(\Omega)$.
Thus we have completed the proof of Lemma \ref{inhomogeneous Boggio 1}.

\medskip
Now using the lemma given above and Theorem \ref{special:helm}, we have the following {\sl completeness} result for homogeneouse viscoelastic system. Let $\Omega\subset {\mathbb R}^3$ be the cylinder given by 
$\Omega = \{ (r,\theta,z) :\, 0 \leq r \leq a, 0 \leq \theta \leq 2 \pi, 0 \leq z \leq h \}$ 
in terms of cylindrical coordinates $(r,\theta,z)$, where $a,\,h$ are positive constants.

\begin{theorem}\label{theorem:StokesHelmholtz:Decompo:Laplace} If $u \in W_p^6(\Omega)$ is a solution of $L_{\lambda,\mu}u=0$ in $\Omega$, then 
$u = \nabla \varphi_1 + \nabla \times (\psi_1 e_z) + \nabla \times \nabla \times (\chi_1 e_z)$, 
where $e_z$ is the unit vector in the $z$ direction, and $\varphi_1$, $\psi_1$, $\chi_1$ are solutions of 
\begin{equation}\label{StokesHelmholtz:Decompo:Laplace}
\left\{
\begin{array}{lll}
(\Delta + \frac{\rho \omega^2}{\lambda + 2\mu}) \varphi_1 &=& 0, \\
(\Delta + \frac{\rho \omega^2}{\mu}) \psi_1 &=& 0, \\
(\Delta + \frac{\rho \omega^2}{\mu}) \chi_1 &=& 0 
\end{array}
\right.
\end{equation}
in $\Omega$, where $\varphi_1 \in H^5(\Omega)$ and $\psi_1 \in H^{5,2}(\Omega)$, $\chi_1 \in H^{6,3}(\Omega)$.
\end{theorem}
\noindent
{\bf Proof.} Since the solution $u$ is in $W_p^6(\Omega)$, we have $u = \nabla \varphi + \nabla \times (\psi e_z) + \nabla \times \nabla \times (\chi e_z)$ in $\Omega$ with $\varphi \in H_{loc}^7({\Bbb R}^3)$, $\psi \in H_{loc}^{7,4}({\Bbb R}^3)$, $\chi \in H_{loc}^{8,5}({\Bbb R}^3)$ by extending $u$ to $W_p^6({\Bbb R}^3)$ with compact support (see \cite{Renardy} about this extension) and Theorem \ref{special:helm}. Then $L_{\lambda, \mu}u=0 $ in $\Omega$ has the form
\begin{eqnarray}\label{phi,psi,chi} 
(\lambda + 2\mu) (\nabla \square_1 \varphi) + \mu \nabla \times (\square_2 \psi e_z) + \mu \nabla \times \nabla \times (\square_2 \chi e_z) = 0 \text{ in } \Omega
\end{eqnarray}
where $\square_1 = \Delta + \frac{\rho \omega^2}{\lambda + 2\mu}$, $\square_2 = \Delta + \frac{\rho \omega^2}{\mu}$. 
In the rest of the proof we will omit ''in $\Omega$" for simplicity. Taking the divergence of \eqnref{phi,psi,chi}, we have
$(\lambda + 2\mu) (\nabla^2 \square_1 \varphi) = 0$.
Here we decompose  $\varphi$ using Lemma \ref{inhomogeneous Boggio 1} into $\varphi = \varphi_1 + \varphi_2$ 
with $\varphi_1,\, \varphi_2 \in H^5(\Omega)$ such that $\square_1 \varphi_1 = 0$, $\nabla^2 \varphi_2 = 0$.
\par

On the other hand, taking curl of \eqnref{phi,psi,chi}, we have $\mu \nabla \times \nabla \times(\square_2 \psi e_z + \nabla \times (\square_2 \chi e_z) = 0$. 
Taking one more curl there, we have $-\mu \nabla \times \nabla^2 (\square_2 \psi e_z + \nabla \times (\square_2 \chi e_z)) = 0$ 
from the identity $\nabla \times \nabla \times \nabla \times = - \nabla \times \nabla^2$.  
Let $\Psi:= \nabla \times  (\psi e_z) + \nabla \times \nabla \times  (\chi e_z)$. Then $\Psi\in H^{6,4}(\Omega)$ and $\nabla^2\square_2\Psi=0$. 
Now we decompose  $\Psi$ using Lemma \ref{inhomogeneous Boggio 1} into $\Psi = \Psi_1 + \Psi_2$ with $\Psi_1,\, \Psi_2 \in H^{4,2}(\Omega)$ such that 
\begin{eqnarray}\label{bogtype-1}
\square_2 \Psi_1 = 0,& \nabla^2 \Psi_2 = 0 
\end{eqnarray}
and $\Psi_1=-\mu(\rho\omega^2)^{-1}\nabla^2\Psi$, $\Psi_2=\mu(\rho\omega^2)^{-1}\square_2\Psi$. That is $\Psi_1,\,\Psi_2$ are given by
\begin{eqnarray*}
\Psi_1&=&-\mu(\rho\omega^2)^{-1}\nabla^2\left(\nabla \times  (\psi e_z) + \nabla \times \nabla \times  (\chi e_z)\right)\\
&=&-\mu(\rho\omega^2)^{-1}\left[ \nabla\times (\nabla^2\psi e_z) + \nabla \times \nabla \times  (\nabla^2  \chi e_z)\right],\\
\Psi_2&=&\mu(\rho\omega^2)^{-1}\square_2\left(\nabla \times  (\psi e_z) + \nabla \times \nabla \times  (\chi e_z)\right)\\
&=&-\mu(\rho\omega^2)^{-1}\left[ \nabla\times (\square_2\psi e_z) + \nabla \times \nabla \times  (\square_2\chi e_z)\right].
\end{eqnarray*}
Hence by defining
\begin{eqnarray}
\psi_1:=-\mu(\rho\omega^2)^{-1} \nabla^2\psi \in H_{loc}^{5,2}({\Bbb R}^3), & \chi_1:=-\mu(\rho\omega^2)^{-1}\nabla^2\chi\in H_{loc}^{6,3}({\Bbb R}^3), \label{defphch-1}\\
\psi_2:=-\mu(\rho\omega^2)^{-1} \square_2\psi\in H_{loc}^{5,2}({\Bbb R}^3) , & \chi_2:=-\mu(\rho\omega^2)^{-1}\square_2\chi\in H_{loc}^{6,3}({\Bbb R}^3),\label{defphch-2}
\end{eqnarray}
we clearly have
\begin{eqnarray}
\nabla\times\left(\square_2\psi_1 e_z+\nabla\times(\square_2\chi_1 e_z)\right)&=&0,\label{psi,chi,1-1}\\
\nabla\times\nabla^2\left(\psi_2 e_z+\nabla\times(\chi_2 e_z)\right)&=&0.\label{psi,chi,2-1}
\end{eqnarray}

Put $\Xi:=\psi_1 e_z+\nabla\times(\chi_1 e_z) \in H_{loc}^{5,2}({\Bbb R}^3)$. Then from \eqref{psi,chi,1-1}, 
we have $\nabla \times \square_2 \Xi=0$ and hence $\square_2 \Xi=\nabla{h}$ for some scalar potential $h\in H^{4,1}({\Bbb R}^3)$. Here we consider the solvability for the equation $\square_2 \Xi_2=\nabla{h} \in H_{loc}^{3,0}({\Bbb R}^3)$. We can show that $\Xi_2:=\square_2^{-1}(\nabla h) \in H^{5,2}({\Bbb R}^3)$ (see \cite{yafaev}). Define $\Xi_1:=\Xi-\Xi_2\in H^{5,2}(\mathbb{R}^3)$. Then, we can have a Boggio type decomposition for $\Xi$ as
$\Xi=\Xi_1+\Xi_2$ with $\square_2 \Xi_1=0$ and $\nabla\times \Xi_2=0$.

Taking curl on $\Xi_1$, we have
\begin{eqnarray}
\nabla\times \Xi_1 &=&\nabla\times(\Xi - \Xi_2) \\ \nonumber
&=&\nabla\times (\psi_1 e_z)+\nabla\times\nabla\times(\chi_1 e_z)\label{special:helm-Bogg-1}
\end{eqnarray}
with 
\begin{equation}
\label{special:helm-Bogg-3}\left\{
\begin{array}{ccc}
\Delta^\prime \psi_1 =&-e_z\cdot(\nabla\times\nabla\times \Xi_1 ),&\\
\Delta^\prime \chi_1=&-e_z\cdot(\nabla\times \Xi_1).&
\end{array}
\right.
\end{equation}

Taking $\square_2$ on (\ref{special:helm-Bogg-3}), we have 
\begin{equation*}
\left\{
\begin{array}{ccc}
\square_2 \Delta^\prime \psi_1 =&0,&\\
\square_2 \Delta^\prime \chi_1=&0.&
\end{array}
\right.
\end{equation*}
Here, we will use the uniqueness for $\Delta^\prime$ in Theorem \ref{special:helm}. To do so, we have to show $\square_2 \psi_1$, $\square_2 \chi_1 \in M_{2,\delta} ({\Bbb R}^2)$ where $-1< \delta <0$ as in the condition (\ref{condition for indices 2}).
From (\ref{defphch-1}), we have $\square_2 \psi_1 \in H_{loc}^{3,0}({\Bbb R}^3) \subset L^2 ({\Bbb R}^3)$ and $\square_2 \chi_1 \in H_{loc}^{4,1}({\Bbb R}^3) \in L^2 ({\Bbb R}^3)$. Since $\delta<0$, we have $L^2 ({\Bbb R}^3) \subset L_{\delta}^2 ({\Bbb R}^3)$. 
So, we have for fixed $z \in {\Bbb R}$, $\square_2 \psi_1 (x',z), \square_2 \chi_1 (x',z) \in L_{\delta}^2 ({\Bbb R}^2)$. Obviously, we have $\Delta{'} \square_2 \psi_1 = \Delta{'} \square_2 \chi_1 = 0 \in L_{\delta +2}^2 ({\Bbb R}^2) $. By applying the lemma \ref{L,M relation}, for any fixed $z \in {\Bbb R}$, we have $\square_2 \psi_1 (x',z), \square_2 \chi_1 (x',z) \in M_{2,\delta}({\Bbb R}^2)$. So, the uniqueness for $\Delta^{'}$ give us that 
\begin{equation*}
\left\{
\begin{array}{ccc}
\square_2 \psi_1 =&0,&\\
\square_2  \chi_1=&0.&
\end{array}
\right.
\end{equation*} 
\par

Then it follows from the decompositions on $\varphi, \psi, \chi$ and (\ref{phi,psi,chi}) that we have
\begin{eqnarray*} 
0 &=& (\lambda + 2\mu) (\nabla \square_1 \varphi) + \mu \nabla \times (\square_2 \psi e_z) + \mu \nabla \times \nabla \times (\square_2 \chi e_z) \\ \nonumber 
&=& (\lambda + 2\mu) (\nabla \square_1 (\varphi_1 + \varphi_2)) + \mu \nabla \times (\square_2 (\psi_1 + \psi_2) e_z) + \mu \nabla \times \nabla \times (\square_2 (\chi_1 + \chi_2) e_z) \nonumber \\
&=& (\lambda + 2\mu) (\nabla \square_1 \varphi_2) + \mu \nabla \times (\square_2  \psi_2 e_z) + \mu \nabla \times \nabla \times (\square_2  \chi_2 e_z).
\end{eqnarray*}
And using $\nabla^2 \varphi_2=0$ and (\ref{psi,chi,2-1}), we have
$$\rho \omega^2 (\nabla \varphi_2) + \rho \omega^2 \nabla \times ( \psi_2 e_z) + \rho \omega^2 \nabla \times \nabla \times (\chi_2 e_z) = 0.$$
It follows that our solution $u$ has the form
\begin{eqnarray*} u &=& \nabla \varphi + \nabla \times (\psi e_z) + \nabla \times \nabla \times (\chi e_z) \\
&=& \nabla \varphi_1 + \nabla \times (\psi_1 e_z) + \nabla \times \nabla \times (\chi_1 e_z) \\
&\quad& + \nabla \varphi_2 + \nabla \times (\psi_2 e_z) + \nabla \times \nabla \times ( \chi_2 e_z) \\ 
&=& \nabla \varphi_1 + \nabla \times (\psi_1 e_z) + \nabla \times \nabla \times (\chi_1 e_z),
\end{eqnarray*}
where the three scalar functions $\varphi_1$, $\psi_1$ and $\chi_1$ satisfy \eqref{StokesHelmholtz:Decompo:Laplace} in $\Omega$. Thus we have completed the proof of Theorem \ref{theorem:StokesHelmholtz:Decompo:Laplace}.

\section{Solution formula to torsional deformation}$\label{Sec:sol-for-tor}$

In this section we will give the solution for a torsional deformation of cylinder. For that we consider again the viscoelastic equations of system 
\begin{eqnarray}\label{homo-viscoelastic-eqn}
L_{\lambda,\mu} u = 0 \quad \text{in} \quad \Omega, 
\end{eqnarray}
where $\Omega\subset {\mathbb R}^3$ is the cylinder given by $\Omega = \{ (r,\theta,z) :\, 0 \leq r \leq a, 0 \leq \theta \leq 2 \pi, 0 \leq z \leq h \}$.
Decompose the boundary $\partial\Omega$ of $\Omega$ into $D= \{ (r,\Gt) ~|~ \, 0 \leq r < a, \ 0 \leq \Gt \leq 2 \pi \}$, a$D_0=D\times\{0\}$ and $D_h=D\times\{h\}$ which are the bottom and top parts of the boundary $\p\GO$, respectively. Then the mixed type boundary condition for torsional deformation is given as 
\beq\label{traction}
\begin{cases}
t_{rr} =  0  \\
t_{r\theta} = 0 \quad\mbox{on } \p\GO \setminus\overline{ (D_0 \cup D_h)} \mbox{,} \\
t_{rz} = 0
\end{cases}
\begin{cases}
t_{zr} = 0  \\
t_{z\theta} = f(r)  \quad\mbox{on } D_0, \\
t_{zz} =  0
\end{cases} \\
\Bu(r,\theta,z) = 0 \quad\mbox{on } D_h.
\eeq
Here we used the notations $t_{ij}$'s to denote tractions in cylindrical coordinates whose definitions can be found in \cite[p. 74-75]{J.D.Achenbach} (see the the next page for precise definitions). The function $f$ represents the torque given at the bottom surface of the cylinder. We assume that it is smooth enough and depends only on the radial variable.
\par
The aim of this section is to give the form of solution $u$ of \eqref{homo-viscoelastic-eqn} satisfying the mixed type boundary type boundary condition \eqref{traction} by using Theorem \ref{theorem:StokesHelmholtz:Decompo:Laplace}. To be precise, we can see from \eqnref{curlsol} given later that the solution $u$ has a form of 
$$u= \nabla \times (\psi(r,z) e_z) $$
for some scalar function $\psi(r,z)$ satisfying \eqnref{StokesHelmholtz:Decompo:Laplace}. This type is well known as torsional wave which involves a displacement in the circumferential direction only (see \cite{Billingham}).
\par
We recall the definitions of $t_{ij}$'s.
For a vector valued function $u=(u_1,u_2, u_3)$, we define $u_r$, $u_\Gt$ and $u_z$ by
\beq\label{urutheta}
u_r := u_1 \cos \Gt + u_2 \sin \Gt , \quad u_\Gt := -u_1 \sin \Gt + u_2 \cos \Gt , \quad u_z :=u_3.
\eeq
Then the components of the traction, denoted by $t_{ij}$ in \eqnref{traction}, are given by the following formulae:
\begin{align} 
t_{rr} &= \Gl \nabla \cdot u + 2 \mu \p_r u_r, \nonumber \\
t_{r\Gt} &= \mu \left( \p_r u_\Gt - r^{-1} u_\Gt + r^{-1} \p_\Gt u_r \right) , \nonumber \\
t_{rz}&= t_{zr} = \mu \left( \p_z u_r + \p_r u_z \right), \label{tdef} \\
t_{z\Gt} &= \mu \left( \p_z u_\Gt + r^{-1} \p_\Gt u_z \right), \nonumber \\
t_{zz} &= \Gl \nabla \cdot u + 2 \mu \p_z u_z. \nonumber
\end{align}
\par
Now we define a basis function for series solution of \eqnref{viscosystem} satisfying \eqnref{traction}. 
Let $J_1$ be the Bessel function of order 1, {\it i.e.}, the solution to
\beq
r^2 J_1''(r) + r J_1'(r) + (r^2 - 1) J_1(r) = 0,
\eeq
and let $0 <k_2 <k_3 < \cdots$ be positive numbers such that
\beq\label{knzero}
k_n J_1' (k_n) = J_1 (k_n).
\eeq
Define $\Gvf_n$ by
\beq\label{c_n}
\Gvf_n(r) =
\begin{cases}
c_1 r, \quad    &n=1, \\
c_n J_1(k_n r), \quad  &n \geq 2,
\end{cases}
\eeq
where $c_n$ is a normalization constant such that
\beq
\| \Gvf_n \|^2 := \int_{0}^{1} |\Gvf_n(r)|^2 r dr =1.
\eeq
We note that $\Gvf_n$ ($n \ge 2$) satisfies
\beq\label{Gvfn4}
r^2 \Gvf_n''(r) +  r \Gvf_n'(r) + (k_n^2 r^2 - 1) \Gvf_n(r) = 0,
\eeq
and
\beq\label{Gvfn2}
\Gvf_n'(1) = \Gvf_n(1).
\eeq
It is known that $\{ \Gvf_n \}_{n=1}^{\infty}$ forms a complete orthonormal system in $L^2((0,1), rdr)$ (see \cite{watson}). 
\par
Let $a$ be the radius of the cylinder, and define
\beq
\Gvf_n^a(r):= a^{-1} \Gvf_n (a^{-1} r).
\eeq
Then $\{ \Gvf_n^a \}_{n=1}^{\infty}$ is a complete orthonormal system in $L^2((0,a), rdr)$. So, the torque $f$ can be expanded as
\beq\label{dini series}
f(r) = \sum_{n=1}^{\infty} f_n \Gvf_n^a(r).
\eeq
Moreover, one can see from \eqnref{Gvfn4} that the following holds:
\beq\label{Gvfna}
r^2 (\Gvf_n^a)''(r) +  r (\Gvf_n^a)'(r) + \left( \left( \frac{k_n}{a} \right)^2 r^2 - 1 \right) \Gvf_n^a(r) = 0.
\eeq

Let $k^2 = \Gr \Go^2 /\mu$ and define $\Gg_n$ by
\beq\label{Ggn}
\Gg_n^2 = \begin{cases}
-k^2, \quad    &n=1, \\
k_n^2/a^2- k^2, \quad  &n \geq 2.
\end{cases}
\eeq
Also let
\beq\label{qnz}
q_n(z) = \frac{e^{\Gg_n(h-z)} - e^{\Gg_n(z-h)}}{ e^{-\Gg_n h  } + e^{\Gg_n h}}, \quad n=1,2,\cdots .
\eeq
Then we have the following lemma and theorem (see \cite[Section2]{EKNW} for their proofs). 

\begin{lemma}\label{lem:soln}
For $n=1, 2, \cdots$ let
\beq\label{vrz}
v_n(r,z) = \frac{1}{\Gm\Gg_n}  q_n(z) \Gvf_n^a(r)
\eeq
and
\beq
\GT(\Gt):=\sin \Gt \, \Be_1 - \cos \Gt \, \Be_2,
\eeq
where $\Be_1=(1,0,0)$, $\Be_2=(0,1,0)$. 
Then $w_n (r, \Gt, z):= v_n(r,z) \GT(\Gt)$  is the solution of mixed type boundary value problem \eqnref{homo-viscoelastic-eqn}, \eqnref{traction} with $f$ replace by $f_n\Gvf_n^a(r)$.
\end{lemma}

\begin{theorem}\label{thm:main}
Let $f=\sum_{n=1}^{\infty} f_n \Gvf_n^a(r)$ satisfy 
\beq\label{H1/2}
\sum_{n=1}^\infty \frac{|f_n|^2}{n} < \infty,
\eeq
or equivalently $f (r, \Gt):=f(r) \GT(\Gt)$ belong to $H^{-1/2}(D_0)$. Then $u$ defined by
\beq\label{solexpress}
u = \sum_{n=1}^\infty f_n w_n =\sum_{n=1}^\infty f_n v_n(r,z) \GT(\Gt)
\eeq
is the unique solution of the mixed type boundary value problem \eqnref{homo-viscoelastic-eqn}, \eqnref{traction} in $H^1(\GO)$. Further $u$ satisfies the estimate
\beq\label{regularity}
\| u \|_{H^1(\GO)} \le C \| f \|_{H^{-1/2}(D_0)}.
\eeq
with some constant $C$ independent of $f$.
\end{theorem}

We remark on Theorem \ref{thm:main} that if we define $\Gy$ by
\beq \nonumber
\Gy (r,z) := \sum_{n=1}^{\infty} \left( \int_0^r \Gvf_n^a(s) ds \right) \frac{f_n}{\Gm \Gg_n}q_n(z),
\eeq
then $u$ can be expressed as
\beq\label{curlsol}
u = \nabla \times (\Gy e_z).
\eeq
where $\Gy$ satifies \eqnref{StokesHelmholtz:Decompo:Laplace}. This is really giving the formula of solution based on the completeness of representation.

\vskip0.7cm
\noindent
{\bf\large Acknowledgement} 
\newline
The second author was partially supported by grant-in-aid for Scientific Research (15K21766 and 15H05740) of the Japan Society for the Promotion of Science doing the research of this paper.

\end{document}